\documentstyle{article}

\newcommand{\be}{\begin{equation}}
\newcommand{\ee}{\end{equation}}
\newcommand{\ba}{\begin{eqnarray}}
\newcommand{\ea}{\end{eqnarray}}
\newcommand{\baa}{\begin{eqnarray*}}
\newcommand{\eaa}{\end{eqnarray*}}
\newcommand{\bb}{\begin{thebibliography}}
\newcommand{\eb}{\end{thebibliography}}

\newcommand{\bi}[1]{\bibitem{#1}}
\newcommand{\lab}[1]{\label{#1}}
\newcommand{\re}[1]{(\ref{#1})}

\newcounter{my}
\newcommand{\he}%
   {\stepcounter{equation}\setcounter{my}%
   {\value{equation}}\setcounter{equation}0%
   }%
\newcommand{\she}%
   {\setcounter{equation}{\value{my}}%
    }%

\newcommand\cn{\mbox{cn}}
\newcommand\dn{\mbox{dn}}
\newcommand\sn{\mbox{sn}}

\newtheorem{th}{Theorem}
\newtheorem{pr}{Proposition}

\newtheorem{lem}{Lemma}

\begin{document}

%\begin{titlepage}
\vspace*{10mm}

\begin{center}

{\Large \bf Elliptic polynomials orthogonal on the unit circle
with a dense point spectrum}

\vspace{5mm}

{\large \bf Alexei Zhedanov}

\medskip

{\em Donetsk Institute for Physics and Technology, Donetsk 83114,  Ukraine}

\end{center}

%\vspace*{5mm}

\begin{abstract}
We introduce two explicit examples of polynomials orthogonal on the
unit circle. Moments and the reflection coefficients are expressed
in terms of Jacobi elliptic functions. We find explicit expression
for these polynomials in terms of a new type of elliptic
hypergeometric function. We show that obtained polynomials are
orthogonal on the unit circle with respect to a dense point meausure, i.e. the spectrum 
consists from infinite number points of increase which are dense on the unit circle. We
construct also corresponding explicit systems of
polynomials orthogonal on the interval of the real axis with respect to a dense point measure. 
They can be considered as an elliptic generalization of the Askey-Wilson polynomials of a special type. 

\end{abstract}

\bigskip\bigskip

%\end{titlepage}

\newpage
\section{Introduction}
Let $\cal L$ be some linear functional defined on all possible
monomials $z^n$ by the moments \be c_n = {\cal L} \{z^n\}, \quad n
=0 , \pm 1, \pm 2 \dots \lab{def_mom} \ee (in general the moments
$c_n$ are arbitrary complex numbers). Then the functional $\cal L$
is defined on the space of generic Laurent polynomials ${\cal
P}(z)= \sum_{n=-N_1}^{N_2}a_n z^n$ where $a_n$ are arbitrary
complex numbers and $N_{1,2}$ are arbitrary integers. Namely we
have (by linearity of the functional) $${\cal L}\{{\cal P}(z)\}=
\sum_{n=-N_1}^{N_2}a_nc_n.$$

The monic Laurent biorthogonal polynomials $P_n(z)$ are defined by
the determinantal formula \cite{HR} \ba \nonumber
&&P_n(z)=(\Delta_n)^{-1} \left |
\begin{array}{cccc} c_0 & c_1 & \dots & c_n \\ c_{-1} & c_0 &
\dots & c_{n-1} \\ \dots & \dots & \dots & \dots\\ c_{1-n}&
c_{2-n}& \dots & c_1\\ 1& z & \dots & z^n  \end{array} \right |,
\lab{deterP} \ea where $\Delta_n$ are defined as \be \Delta_n =
\left | \begin{array}{cccc} c_0 & c_{1} & \dots &
c_{n-1}\\ c_{-1}& c_0 & \dots & c_{n-2}\\ \dots & \dots & \dots & \dots\\
c_{1-n} & c_{2-n} & \dots & c_0 \end{array} \right | \lab{Delta}
\ee It is obvious from the definition \re{deterP} that the
polynomials $P_n(z)$ satisfy the orthogonality property \be {\cal
L}\{P_n(z) z^{-k}\} =h_n \delta_{kn}, \quad 0 \le k \le n,
\lab{ort_P} \ee where the normalization constants are \be h_0=c_0,
\quad h_n = \Delta_{n+1}/\Delta_{n}. \lab{h_n} \ee This
orthogonality property can be rewritten in terms of biorthogonal
relation \cite{Pas}, \cite{HR}, \be {\cal L}\{P_n(z) Q_m(1/z)\} =
h_n \delta_{nm} \lab{biPQ} \ee where $h_n = \Delta_{n+1}/\Delta_n$
and the polynomials $Q_n(z)$ are defined by the formula \ba
&&Q_n(z)=(\Delta_n)^{-1} \left |
\begin{array}{cccc} c_0 & c_{-1} & \dots & c_{-n} \\ c_{1} & c_0 &
\dots & c_{1-n} \\ \dots & \dots & \dots & \dots\\ c_{n-1}&
c_{n-2}& \dots & c_{-1}\\ 1& z & \dots & z^n  \end{array} \right
|, \lab{deterQ} \ea i.e. the polynomials $Q_n(z)$ are again
Laurent biorthogonal polynomials  with moments $c^{\{Q\}}_n =
c_{-n}$.

Laurent biorthogonal polynomials satisfy three-term recurrence
relation \cite{HR} \be P_{n+1}(z) + d_n P_n(z) = z(P_n(z) + b_n \:
P_{n-1}(z)) \lab{3LBP} \ee with some recurrence coefficients $b_n,
d_n$. In fact, recurrence relation \re{3LBP} uniquely determine
Laurent biorthogonal polynomials $P_n(z)$ under the standard
initial conditions $P_{-1}=0, \; P_0=1$.

An important special case of the Laurent biorthogonal polynomials
is obtained if \be c_{-n} = \bar c_n \lab{Sz_c} \ee where $\bar
c_n$ means complex conjugated moments with respect to $c_n$.

In this case the biorthogonal partners are $Q_n(z)= \bar P_n(z)$.
If additionally we demand that $$ \Delta_n>0, \quad n=0,1,2,\dots
$$ then there exists a nondecreasing measure $\mu$ on the unit
circle $|z|=1$ such that \be \int_{0}^{2 \pi} P_n(e^{i \theta})
\bar P^m(e^{-i \theta}) d \mu(\theta) = h_n \: \delta_{nm}
\lab{Sz_ort} \ee Corresponding polynomials are called the
Szeg\H{o} polynomials orthogonal on the unit circle. In what
follows we will denote the Szeg\H{o} polynomials by $\Phi_n(z)$
instead of $P_n(z)$.

The Szeg\H{o} polynomials satisfy the recurrence relation \be
\Phi_{n+1}(z) = z\Phi_n(z) - \bar a_{n} z^n \bar \Phi_n(1/z),
\lab{Sz_rec} \ee where $a_n=-\bar \Phi_{n+1}(0)$ are so-called
reflection parameters (sometimes called Schur, Geronimus,
Verblunsky... parameters). The normalization constants $h_n$ in
the orthogonality relation \re{Sz_ort} are expressed in terms of
the reflection parameters as \be  h_n = (1-|a_0|^2) (1-|a_1|^2)
\dots (1-|a_{n-1}|^2) \lab{h_a} \ee If $|a_n| \ne 1$ then the
orthogonal polynomials are nondegenerate, i.e. $h_n \ne 0, \;
n=0,1,\dots$. If however $|a_N|=1$ for some $N>1$ then we have a
{\it finite} system of orthogonal polynomials on the unit circle:
$\Phi_0(z), \Phi_1(z), \dots \Phi_N(z)$.

From the Szeg\H{o} recurrence relation \re{Sz_rec} one can derive
three term recurrence relation \be \Phi_{n+1}(z) + d_n \Phi_{n}(z)
= z(\Phi_{n}(z) + b_n \Phi_{n-1}(z)), \lab{3UC} \ee where
$$
d_n = - \frac{\bar a_n}{\bar a_{n-1}}, \quad b_n = d_n
(1-|a_{n-1}|^2)
$$
Note that relation \re{3UC} is a special case of three-term
recurrence relation \re{3LBP} for the Laurent biorthogonal
polynomials.

In what follows we restrict ourselves with the so-called real
case, i.e. we will assume that $c_0=1$ (the standard normalization
condition) and all the moments $c_n$ are real and symmetric
$c_{-n}=c_n$, whereas the reflection parameters $a_n$ are real and
satisfy the restriction \be -1 < a_n < 1, \quad n=0,1,\dots
\lab{res_a} \ee In this case the polynomials $\Phi_n(z)$ have the
orthogonality property \be \int_{0}^{2 \pi} \Phi_n(e^{i \theta})
\Phi_m(e^{-i \theta}) d \mu(\theta) = h_n \: \delta_{nm},
\lab{ort_real} \ee where the measure possesses symmetric property
$$
\mu(2\pi - \theta) =- \mu(\theta)
$$
In case if the measure admits existence of the weight function
$\rho(\theta)$
$$
d\mu(\theta) = \rho(\theta) d\theta
$$
we have that $\rho(\theta)$ is a real nonnegative function on the
unit circle (i.e. $\rho(\theta)\ge 0, \quad 0\le \theta \le 2\pi$)
which is symmetric with respect to the reflection
$$
\rho(2\pi - \theta ) = \rho(\theta)
$$
Note that in the case \re{res_a} all expansion coefficients of the
polynomials $P_n(z)$ are real, i.e. polynomials $P_n(z)$ take real
values for real $z$.

Orthogonality condition \re{ort_real} in this case is equivalent
to the conditions \be \int_{0}^{2 \pi} \Phi_n(e^{i \theta}) e^{-i
m \theta} d \mu(\theta) =0, \quad m=0,1,\dots n-1 \lab{ort_alt}
\ee The recurrence relation for this case becomes \be
\Phi_{n+1}(z) = z\Phi_n(z) - a_{n} z^n \Phi_n(1/z), \lab{Sz_real}
\ee with the initial condition $\Phi_0(z) =1$. Note that if the
sequence of the real reflection parameters $a_n, \; n=0,1,2,\dots$
is given then the monic polynomials $\Phi_n(z)$ are uniquely
determined through the recurrence relation \re{Sz_real}. Vice
versa, if the real moments $c_n$ are given ($c_0=1$) with the
symmetry condition $c_{-n}=c_n$ and with the positivity condition
$\Delta_n > 0, \; n=1,2,\dots$ then all the reflection
coefficients $a_n$ are determined uniquely and satisfy the
restriction $-1 < a_n <1, \; n=0,1,\dots$.

In what follows we will use a simple
\begin{pr} \lab{pr_c}
Assume that $-1 < a_n <1$ for all $n=0,2,\dots$ and corresponding
polynomials $\Phi_n(z)$ have the expression
$$
\Phi_n(z) = \sum_{s=0}^n W_{ns} z^s
$$
with some expansion coefficients $W_{ns}$ ($W_{nn}=1$ due to
monicity of polynomials $\Phi_n(z)$). Assume that there exists a
real sequence $\tilde c_n, \; n=0,1,2,\dots$ such that $\tilde c_0
=1$ and  \be \sum_{s=0}^n W_{ns} \tilde c_s =0, \quad
n=1,2,3,\dots \lab{cond_Wc} \ee Then the sequence $\tilde c_n$
coincides with the (unique) sequence of moments $c_n$
corresponding to the polynomials $\Phi_n(z)$.
\end{pr}
{\it Proof}. Under condition $-1 <a_n<1$ we have that the
polynomials $\Phi_n(z)$ have a nondecreasing measure on the unit
circle with some unique moments $c_n$ ($c_0=1$). The orthogonality
relation for these polynomials can be presented in the form \be
\langle {\cal L}\Phi_n(z), z^{-m} \rangle =0, \quad m=0,1,2,
\dots, n-1 \lab{ort_Pm} \ee or in an equivalent form
\be \sum_{s=0}^n W_{ns} c_{s-m} =0, \quad n=1,2,\dots, \;
m=0,1,\dots, n-1, \lab{ort_mc} \ee where $c_n = \langle {\cal L},
z^n \rangle$ are the unique moments corresponding to the linear
functional ${\cal L}$. In particular, for $m=0$ we have the
relation \be \sum_{s=0}^n W_{ns} c_{s} =0, \quad n=1,2,\dots,
\lab{ort_m=0} \ee Now consider relation \re{cond_Wc} for $n=1$. We
have $\tilde c_1 -a_0 \tilde c_0=0$. Hence $\tilde c_1 = a_0 =
c_1$ because $\tilde c_0 =1$. We see that the moments $\tilde c_1$
and $c_1$ coincide. In order to prove that $\tilde c_n = c_n$ for
all $n=1,2,3, \dots$ it is sufficient to apply induction with
respect to condition \re{cond_Wc}. This finishes proving of the
proposition. The importance of this proposition is that if we know
the only orthogonality condition \re{cond_Wc} then the general
orthogonality relations \re{ort_mc} and \re{ort_Pm} holds
automatically and moreover the sequence $\tilde c_n$ coincides
with the moment sequence $c_n$.

Together with given polynomials $\Phi_n(z)$ one can introduce the
"reflected sign" polynomials $\tilde \Phi_n(z) = (-1)^n
\Phi_n(-z)$.
\begin{pr}
If $\Phi_n(z)$ are polynomials orthogonal on the unit circle with
real reflection coefficients $a_n$ then the polynomials $\tilde
\Phi_n(z)$ are again polynomials orthogonal on the unit circle
with the reflection coefficients $\tilde a_n = (-1)^{n+1} a_n$.
The expansion coefficients for the polynomials $\tilde \Phi_n(z)$
are $\tilde W_{ns} = (-1)^{n+s} W_{ns}$. The moments $\tilde c_n$
corresponding to the polynomials $\tilde \Phi_n(z)$ are $\tilde
c_n =(-1)^n c_n$.
\end{pr}
The proof of this proposition is trivial. We will use this
proposition in order to reduce the expansion coefficients $W_{ns}$
and reflection parameters $a_n$ to the most convenient form.

Note also that for the case of real parameters $a_n$ we have
simple expressions for $\Phi_n(\pm 1)$: \be \Phi_n(1) =
\prod_{s=0}^{n-1}(1-a_s) \lab{Phi1} \ee and \be \Phi_n(-1) =
(-1)^n \: \prod_{s=0}^{n-1}(1+ (-1)^s a_s) \lab{Phi-1} \ee These
formulas follow directly from recurrence relations \re{Sz_rec}.

\section{Elliptic binomial coefficients and their properties}
\setcounter{equation}{0} Define the so-called "elliptic binomial
coefficients" (EBC) $E^n_j$. (In what follows we will often omit
dependence of the Jacobi elliptic functions like $\sn(z;k)$ on the modulus $k$).
If $0 \le j \le n$ then \be E_0^n = 1, \quad E_j^n =
\prod_{s=0}^{j-1} \frac{\sn(w(n-s))}{\sn(w(s+1))}, \lab{E_def} \ee
and $E^n_j=0$ otherwise. In \re{E_def} $w$ is an arbitrary real
parameter satisfying the restriction \be w \ne 4K M_1/M_2
\lab{restr_w} \ee with some nonzero integers $M_1,M_2$. Recall
that $K \equiv K(k)$ is the complete elliptic integral of the
first kind.

The elliptic modulus $k$ is assumed to belong to the standard
interval $0<k<1$. In this case all EBC are real for all
$j,n=0,1,\dots$.

We have the obvious symmetry property \be E_{n-j}^n = E_j^n, \quad
j=0,1,\dots, n \lab{sym_E} \ee Moreover, in the limit $w=0$ we
obtain the ordinary binomial coefficients
$$
\lim_{w\to 0} E_j^n = {n \choose j} = \frac{n! }{j! (n-j)!}
$$
In the limit $k=0$ we get trigonometric degeneration of the
elliptic Jacobi functions: \be \lim_{k \to 0} E^n_j =
\prod_{s=0}^{j-1} \frac{\sin(w(n-s))}{\sin(w(s+1))}=
q^{\frac{nj-j^2}{2}} \: \left [ {n \atop j}  \right ]_q, \quad
q=e^{-2iw}, \lab{E_k=0} \ee where $\left [ {n \atop j} \right ]_q$
are so-called the Gauss polynomials or q-binomial coefficients
\cite{And}, \cite{KS} defined as \be \left [ {n \atop j}  \right
]_q = \frac{(1-q^n)(1-q^{n-1}) \dots (1-q^{n-j+1})}{(1-q) (1-q^2)
\dots (1-q^j)}= \frac{(q)_n}{(q)_j(q)_{n-j}}, \lab{qGauss} \ee
where $(a)_n =(1-a)(1-aq) \dots (1-aq^{n-1})$ is standard
Pochhammer q-symbol (q-shifted factorial).

Another limit $k=1$ leads to hyperbolic functions degeneration. In
this limit we have \be \lim_{k \to 1} E^n_j =
\frac{(q)_n}{(q)_j(q)_{n-j}} \frac{(-q)_j(-q)_{n-j}}{(-q)_n} =
(-1)^j \frac{(-q)_j (q^{-n})_j}{(q)_j (-q^{-n})_j}, \lab{hyp_E}
\ee where $q=e^{-2w}$.

A less trivial property of the EBC is the recurrence relation: \be
\cn(wj) \dn(w(n-j)) E^n_j = \cn(wn) E^{n-1}_{j-1} + \dn(wn)
E^{n-1}_j \lab{rec_E} \ee In order to prove the relation
\re{rec_E} we need addition formulas for the Jacobi elliptic
functions \cite{WW} from which the identity follows
$$
\sn(wn) \dn(w(n-j)) \cn(wj) = \cn(wn) \sn(wj) + \dn(wn)\sn(w(n-j))
$$
which is valid for all (complex) values of $w,n,j$. Using this
identity we immediately arrive at recurrence relation \re{rec_E}.
Note that in the limit $w=0$ relation \re{rec_E} becomes the well
known recurrence relation for the ordinary binomial coefficients:
\be {n \choose j} = {n-1 \choose j-1} + {n-1 \choose j}
\lab{bin_rec} \ee

There are other recurrence relations for EBC of the similar type
\re{rec_E}. For example if one replaces $j \to n-j$ in \re{rec_E}
and apply the symmetry property \re{sym_E} then one obtains \be
\cn(w(n-j)) \dn(wj) E^n_j = \cn(wn) E^{n-1}_{j} + \dn(wn)
E^{n-1}_{j-1} \lab{rec_E1} \ee which differs from \re{rec_E} by
interchanging of the functions $\cn$ and $\dn$.  We present also
the relation \be \cn(w(n-j)) E^n_j = \dn(w(n-j)) E^{n-1}_{j-1} +
\cn(wn)\dn(wj) E^{n-1}_j \lab{rec_E2} \ee and the similar one \be
\cn(wj) E^n_j = \dn(wj) E^{n-1}_{j} + \cn(wn)\dn(w(n-j))
E^{n-1}_{j-1} \lab{rec_E3} \ee obtained from \re{rec_E2} by the
change $j \to n-j$.

Another pair of relations \be \dn(w(n-j)) E^n_j = \cn(w(n-j))
E^{n-1}_{j-1} + \dn(wn)\cn(wj) E^{n-1}_j \lab{rec_E4} \ee and \be
\dn(wj) E^n_j = \cn(wj) E^{n-1}_{j} + \dn(wn)\cn(w(n-j))
E^{n-1}_{j-1} \lab{rec_E5} \ee is formally obtained from
\re{rec_E2} and \re{rec_E3} by interchanging $\cn$ and $\dn$
functions.

All these relations are derived in a similar way using addition
formulas for the elliptic Jacobi functions. Moreover, all these
relations are reduced to the standard one \re{bin_rec} the limit
$w \to 0$.

Note that our definition of the "elliptic binomial coefficients"
differs from the one given, e.g. in \cite{GR}.

\section{\cn-polynomials on the unit circle}
\setcounter{equation}{0} Introduce the polynomials \be
\Phi_n^{(C)}(z)= \sum_{s=0}^n W_{ns} z^s, \lab{P_W} \ee where the
expansion coefficients $W_{ns}$ are defined as follows. If
$n=0,2,4,\dots$ is even then \be W_{ns}= (-1)^{s} \: \dn(w(n-s))
\: E_s^n \lab{exp_cn_even} \ee If $n=1,3,5,\dots$ is odd then \be
W_{ns} = (-1)^{s+1} \: \cn(w(n-s)) \: E_s^n \lab{exp_cn_odd} \ee
It is seen that $W_{nn}=1$ for all $n=1,2,3,\dots$ and hence the
polynomials $P_n(z) = z^n + O(z^{n-1})$ are monic.

Define the sequence $a_n = - \Phi_{n+1}^{(C)}(0), \; n=0,1,\dots$.
From definition \re{P_W} it follows that \be a_{2n} =
\cn(w(2n+1)), \quad a_{2n+1}= - \dn(w(2n+2)), \quad n=0,1,2,\dots
\lab{a_cn} \ee

For the reciprocal polynomials $\Phi_n^*(z) = z^n \Phi_n(1/z)$ we
have (by symmetry property of the elliptic binomial coefficients)
similar expression
$$
{\Phi_n^{(C)}}^*(z) = \sum_{s=0}^n (-1)^{s} \dn(ws) E_s^n z^s
$$
for even values of $n$ and
$$
{\Phi_n^{(C)}}^*(z) = \sum_{s=0}^n (-1)^{s} \cn(ws) E_s^n z^s
$$
for odd values of $n$.

\begin{pr} The polynomials $\Phi_n^{(C)}(z)$ satisfy the recurrence relation \be \Phi_{n+1}^{(C)}(z) =z\Phi_n^{(C)}(z) - a_n
{\Phi^{(C)}}^*_n(z), \quad n=0,1,2,\dots \lab{rec_PC} \ee \end{pr}

Indeed, comparison of terms $z^s$ leads to the conclusion that for
the  even values of $n$ the recurrence relation \re{rec_PC} is
equivalent to
$$
- \cn(w(n+1-s)) E_s^{n+1} +  \dn(w(n+1-s)) E_{s-1}^n =- a_n
\dn(ws) E_s^n
$$
where $a_n = \cn(w(n+1))$. But this relation coincides with
\re{rec_E2}. Quite analogously, for the odd values of $n$ we see
that condition of vanishing of all terms $z^s$ in expression
\re{rec_PC} is equivalent to relation \re{rec_E4}.

From their definition it follows that the parameters $a_n$ satisfy
an obvious restriction $-1 < a_n <1$ for all $n=0,1,\dots$. Thus
the polynomials $P_n^{(C)}(z)$ belong to the class of the
polynomials orthogonal on the unit circle. From general properties
\cite{Ger} it follows that there exists a positive measure on the
unit circle providing orthogonality property.

We define a linear functional $\sigma$ by its moments \be c_n
\equiv \langle \sigma, z^n \rangle = \cn(wn) \lab{c_cn} \ee

\begin{pr} If the functional $\sigma$ is defined as \re{c_cn} then
the polynomials $\Phi_n^{(C)}(z)$ satisfy the orthogonality
property \be \langle \sigma, \Phi_n^{(C)}(z) \rangle \equiv
\sum_{s=0}^n W_{ns} c_s =0, \quad n=1,2,\dots \lab{ort_c_cn} \ee
\end{pr} {\it Proof}. If $n$ is odd then the proof is trivial.
Indeed, we have
$$S\equiv \sum_{s=0}^n W_{ns}c_s = -\sum_{s=0}^n (-1)^s
\cn(w(n-s)) \cn(ws) E^n_s.$$ Changing summation variable $s \to
n-s$ and using symmetry property \re{sym_E} we have
$$S=-\sum_{s=0}^n (-1)^{n-s} \cn(ws) \cn(w(n-s)) E^n_s = -S,$$
whence $S=0$.

Now consider the case when $n$ is even. In this case we have
$$
S\equiv \sum_{s=0}^n W_{ns}c_s = \sum_{s=0}^n (-1)^s \dn(w(n-s))
\cn(ws) E^n_s
$$
We can apply relation \re{rec_E} to present the sum $S$ in an
equivalent form
$$
S=\cn(wn) \sum_{s=0}^{n} (-1)^s E^{n-1}_{s-1} + \dn(wn)
\sum_{s=0}^{n-1} (-1)^s E^{n-1}_s
$$
But $E^{n-1}_{-1}=0$, hence we can replace $s \to s+1$ in the
first sum, then we have
$$
S= (\dn(wn) - \cn(wn)) \: \sum_{s=0}^{n-1} (-1)^s E^{n-1}_s =0
$$
due to symmetry property \re{sym_E} (note that in this case $n-1$
is the odd and all terms in the sum vanish pairwise). Thus we
proved the proposition for all values $n=1,2,\dots$.

Now apply Proposition \re{pr_c} in order to conclude that
condition \re{ort_c_cn} is equivalent to more general
orthogonality condition \be \langle \sigma, \Phi_n^{(C)}(z) z^{-j}
\rangle =0, \quad j=0,1,\dots, n-1 \lab{ort_j} \ee or,
equivalently to condition \be \langle \sigma, \Phi_n^{(C)}(z)
\Phi_m^{(C)}(1/z) \rangle =h_n \: \delta_{nm}, \lab{ort_nm} \ee
Moreover, we proved that the moments $c_n$ providing orthogonality
relation for the polynomials $P_n^{(C)}(z)$ are given by \re{c_cn}
and the functional $\sigma$ defined by \re{c_cn} coincides with
the orthogonality functional ${\cal L}$.

In order to get the orthogonality measure $\mu(\theta)$ on the
unit circle corresponding to the moments  $c_n = \cn(wn)$ we use
the Fourier expansion of the function $\cn(z)$ \cite{WW}: \be
\cn(z;k) = \frac{\pi}{k K} \: \sum_{s=-\infty}^{\infty}
\frac{\exp(i\pi z(s- 1/2)/K)}{q^{s-1/2} + q^{-s+ 1/2}},
\lab{Four_cn} \ee where \be q = \exp(-\pi K'/K) \lab{qKK} \ee and
$K' \equiv K(k'), \; k'^2 = 1-k^2$.

We have
\begin{th}
The polynomials $\Phi_n(z)$ corresponding to the reflection
parameters \re{a_cn} satisfy the orthogonality property on the
unit circle \be \sum_{s=-\infty}^{\infty} \rho_s
\Phi_n^{(C)}(z_s)\Phi_m^{(C)}(1/z_s) = h_n \delta_{nm},
\lab{ort_cn} \ee where \be z_s = \exp(i \pi w (s-1/2)/K)
\lab{zs_cn} \ee and the discrete weights $\rho_s$ are given by \be
\rho_s = \frac{\pi}{k K} \: \frac{1}{q^{s-1/2} + q^{-s+1/2}}
\lab{W_cn} \ee
\end{th}

{\it Proof}. It is sufficient to show that the measure given by
\re{ort_cn} reproduces "true" moments $c_n$, i.e. we should verify
that \be c_n \equiv \cn(wn,k) = \sum_{s=-\infty}^{\infty} \rho_s
z_s^n, \quad n=0, \pm 1, \pm 2, \dots \lab{vermom_cn} \ee But
property \re{vermom_cn} follows immediately from formula
\re{Four_cn} if we put $z=wn$. Thus we verified that our measure
gives correct values of all moments $c_n = \cn(wn)$. This finishes
the proof, because in our case this measure is the unique due to
the condition $-1 < a_n <1$.

In order to get expression for the normalization constants $h_n$
we use formulas \re{h_a} and \re{a_cn} which yields \be h_n =
\mu_n \prod_{s=1}^n \sn^2(ws), \lab{h_n_cn} \ee where $\mu_n =
k^n$ for the even $n$ and $\mu_n = k^{n-1}$ for the odd $n$.
Clearly, the normalization coefficients are all positive $h_n
>0, \; n=0,1,2,\dots$.

Clearly, for all $s=0, \pm 1, \pm 2, \dots$ the spectral points
(i.e. the only points of growth of the orthogonality measure)
$z_s$ lie on the unit circle. Note that under the condition
\re{restr_w} we have that  $z_s$ are {\it dense} on the unit
circle. This means that for every point $z_s$ and for any small
parameter $\epsilon$ there exists another spectral point $z_{s'}$
such that $|z_s-z_{s'}|< \epsilon$. In recent monograph by Barry
Simon \cite{Simon} one can find general results concerning
existence of the singular continuous and dense point measures on the
unit circle. Our $\cn$-polynomials $\Phi_n(z)$ provide a very
simple explicit example leading to such measures.

Sometimes the "reflected sign" polynomials $\tilde \Phi_n^{(C)}(z)
= (-1)^n \Phi_n^{(C)}(-z)$ are more convenient. Using Proposition
{\bf 2} we find that polynomials $\tilde \Phi_n^{(C)}(z)$ have the
moments $\tilde c_n = (-1)^n \cn(wn)$ and the reflection
parameters $\tilde a_n = -\cn(w(n+1)$ for the even $n$ and $\tilde
a_n = - \dn(w(n+1))$ for the odd $n$. The expansion coefficients
are $\tilde W_{ns} = \dn(w(n-s)) E^n_s$ for even $n$ and $W_{ns}
=\cn(w(n-s)) E^n_s$ for the odd $n$. Transformation
$\Phi_n^{(C)}(z) \to \tilde \Phi_n^{(C)}(z)$ in this case is
equivalent to the substitution $w \to w + 2K$. Indeed, we have
$\tilde c_n = \cn(wn + 2Kn) = (-1)^n \cn(wn)$ due to the property
$\cn(z+2K) =-\cn(z)$. This observation allows one to construct
corresponding orthogonality relation for the polynomials $\tilde
\Phi_n^{(C)}(z)$. We have almost the same relation as \re{ort_cn}:
\be \sum_{s=-\infty}^{\infty} \rho_s \tilde \Phi_n^{(C)}(\tilde
z_s) \tilde \Phi_m^{(C)}(1/\tilde z_s) = h_n \delta_{nm},
\lab{ort_cn1} \ee with the same weights $\rho_s$ and normalization
coefficients and where \be \tilde z_s = \exp(i \pi w (s-1/2)/K - i
\pi) = - z_s \lab{zs_cn1} \ee

\section{\dn-elliptic polynomials}
\setcounter{equation}{0} In the previous section we considered the
polynomials $\Phi_n^{(C)}(z)$ orthogonal on the unit circle with
respect to the moments sequence $c_n = \cn(wn), n=0, \pm 1, \pm 2,
\dots$. Consider now the polynomials $\Phi_n^{(D)}(z)$
corresponding to the  moments given by the expression \be c_n=
\dn(wn), \quad n=0,\pm 1, \pm 2, \dots \lab{mom_dn} \ee

Clearly, the moments \re{mom_dn} are real and satisfy the symmetry
property $c_{-n}=c_n$ (it is assumed that $0 <k<1$).

As in the previous section it is verified that the reflection
parameters $a_n$ corresponding the moments \re{mom_dn} are given
by the expression \be a_{2n} = \dn(w(2n+1);k), \quad a_{2n+1}= -
\cn(w(2n+2);k), \quad n=0,1,2,\dots \lab{a_dn} \ee

The explicit expression of the polynomials $\Phi_n^{(D)}(z) =
\sum_{s=0}^n W_{ns} z^s$ in this case is very similar to
expression for polynomials $\Phi_n^{(C)}(z)$. If $n=0,2,4,\dots$
is even then \be W_{ns}= (-1)^{s} \: \cn(w(n-s)) \: E_s^n
\lab{exp_dn_even} \ee If $n=1,3,5,\dots$ is odd then \be W_{ns} =
(-1)^{s+1} \: \dn(w(n-s)) \: E_s^n \lab{exp_dn_odd} \ee Hence in
this case the coefficients $W_{ns}$ are obtained from the
coefficients \re{exp_cn_even} and \re{exp_cn_odd} by simple
interchanging the functions $\cn$ and $\dn$.

For real values of the parameter $w$ the reflection parameters
$a_n$ are real and satisfy the property $-1 \le a_n \le 1$. If
additionally the parameter $w$ is chosen such that $ w \ne 2KM/N $
then we have $-1 < a_n < 1$ for all $n=0,1,2,\dots$ and hence the
reflection parameters satisfy the condition \re{res_a} meaning
that there exists a symmetric measure $\mu(\theta)$ providing
orthogonality of the polynomials $\Phi_n^{(D)}(z)$ on the unit
circle.

In a similar way as in the previous section we obtain the
orthogonality measure
\begin{th}
The polynomials $\Phi_n^{(D)}(z)$ corresponding to the reflection
parameters \re{a_dn} satisfy the orthogonality property on the
unit circle \be \sum_{s=-\infty}^{\infty} \rho_s
\Phi_n^{(D)}(z_s)\Phi_m^{(D)}(1/z_s) = h_n \delta_{nm},
\lab{ort_dn} \ee where \be z_s = \exp(i \pi w s/K) \lab{zs_dn} \ee
and the discrete weights $\rho_s$ are given by \be \rho_s =
\frac{\pi}{ K} \: \frac{1}{q^{s} + q^{-s}} \lab{W_dn} \ee with $q$
given by \re{qKK}.
\end{th}
The proof of this theorem follows immediately from the Fourier
expansion of the Jacobi $\dn$-function \cite{WW}
 \be \dn(z;k) = \frac{\pi}{K} \:
\sum_{s=-\infty}^{\infty} \frac{\exp(i\pi zs/K)}{q^{s} + q^{-s}},
\lab{Four_dn} \ee

\section{Elliptic derivative operator}
\setcounter{equation}{0}
The polynomials $\Phi_n^{(C)}(z)$ and
$\Phi_n^{(D)}(z)$ are very close to one another in their explicit
expression. There is a remarkable relation between them based on a
notion of so-called "elliptic derivative operator" ${\cal E}$. We
define the operator ${\cal E}$ on the space of all formal series
$\sum_{s=0}^{\infty} c_s z^s$ by its action on monomials: \be
{\cal E} z^n = e_n z^{n-1} \lab{Ezn} \ee where
$$
e_n \equiv E^n_1 =\frac{\sn(wn)}{\sn(w)}
$$
is the "elliptic number". Note that in the limit $w=0$ we have
$e_n =n$ and the operator ${\cal E}$ becomes the ordinary
derivative operator: ${\cal E}_{w=0} = \partial_z$. In the limit
$k=0$ we get the q-derivative operator
$$
{\cal E} z^n = \frac{q^n - q^{-n}}{q-q^{-1}} z^{n-1}
$$
where $q=\exp(i w)$.

The operator ${\cal E}$ can be presented explicitly as an infinite
sum of q-derivative operators. Indeed, introduce a set of
q-derivative operators ${\cal W}_j$ which act by the formula \be
{\cal W}_j f(z) = \frac{f(z e^{ij \theta}) -  f(z e^{-ij
\theta})}{2iz}, \quad j=1,2,3,\dots \lab{w-deriv} \ee with some
fixed real parameter $\theta$ (it is assumed that $\theta \ne 2
\pi M/N$ with integer $M,N$). On the set of monomials $z^n$ the
operator ${\cal W}_j$ acts as \be {\cal W}_j z^n = \sin(nj \theta)
z^{n-1}, \quad n=0,1,2,\dots \lab{W_mon} \ee Hence the operator
${\cal W}_j$ sends any polynomial of degree $n$ to a polynomial of
degree $n-1$.

Now we have \be {\cal E} = \sum_{j=1}^{\infty} \beta_j {\cal
W}_{2j-1}, \lab{expl_E_W} \ee where $\theta = \pi w/(2 K)$ and
$$
\beta_j = \frac{2 \pi}{K k  (q^{1/2-j} - q^{j-1/2}) \: \sn(w)},
\quad q = \exp(-\pi K'/K)
$$
Indeed, formula \re{expl_E_W} follows from the Fourier expansion
of the elliptic $\sn$ function \cite{WW}: $$\sn(z;k) = \frac{2
\pi}{Kk} \: \sum_{j=1}^{\infty} \frac{\sin((2j-1)t)}{q^{1/2-j} -
q^{j-1/2}}, \quad t= \pi z/(2K)$$

The operator ${\cal E}$ transforms any polynomial in $z$ of degree
$n$ to a polynomial of degree $n-1$.

We have
\begin{pr}
The polynomials $\Phi_n^{(C)}(z)$ and $\Phi_n^{(D)}(z)$ are
related as \be e_n \: \Phi_{n-1}^{(D)}(z)= {\cal E}
\Phi_n^{(C)}(z), \quad e_n \: \Phi_{n-1}^{(C)}(z)= {\cal E}
\Phi_n^{(D)}(z) \lab{EPP} \ee
\end{pr}
The proof is elementary and is based on explicit expression of the
polynomials $\Phi_n^{(C)}(z),\; \Phi_n^{(D)}(z)$.

As an elementary consequence of this proposition we obtain that
the operator ${\cal E}^2$ transforms both families of polynomials
$\Phi_n^{(C)}(z)$ and $\Phi_n^{(D)}(z)$ to themselves:
$$
{\cal E}^2 \Phi_n^{(C)}(z) = e_n e_{n-1} \Phi_{n-2}^{(C)}(z),
\quad {\cal E}^2 \Phi_n^{(D)}(z) = e_n e_{n-1} \Phi_{n-2}^{(D)}(z)
$$

\section{Limits $k=0$ and $k=1$}
\setcounter{equation}{0} The "trigonometric" limit $k =0$ is not
interesting because it leads to a degeneration. Indeed, in this
limit, e.g. for the $\cn$ polynomials we have $a_{2n+1}=-1$ for
all $n$. Similarly, for $\dn$ polynomials in this limit we have
$a_{2n}=1$. This means that we have the degenerated case of the
polynomials orthogonal on the unit circle. It is instructive to
consider what happens with orthogonality relation in these
degeneration cases. In the limit $k=0$ we have \cite{WW} $K(0) =
\pi/2$ and $K' = \log(4/k) + O(k^2 \log(k))$. Hence $q \to k^2/16
\to 0$. For the case of $cn$-polynomials we have that only terms
$\rho_0$ and $\rho_1$ in formula \re{W_cn} remain nonzero in the
limit $k=0$. Indeed, in this limit we have only two spectral
points $z_0=e^{-iw}$ and $z_1= e^{iw}$ having two equal
concentrated masses $\rho_0 = \rho_1=1/2$. For the case of
$dn$-polynomials we have the only spectral point $z_0=1$ with the
concentrated mass $\rho_0=1$. Thus in the trigonometric limit the
problem becomes trivial.

Consider the "hyperbolic" limit $k=1$. In this case the reflection
parameters $a_n$ for both $cn$ and $dn$ polynomials coincide: \be
a_n = \frac{(-1)^n }{\cosh(w(n+1))} \lab{a_n_hyp} \ee The moments
for both polynomials are  \be c_n = 1/\cosh(wn) \lab{hyp_mom} \ee
Hence both polynomials coincide in this limit and we have from
\re{hyp_E} the following expression for the power coefficients \be
W_{ns} = 2\frac{(-1)^{n+s} q^{-n/2+s/2}}{1+q^{-n}}
 \: \frac{(q^{-n})_s (-q)_s}{(q)_s(-q^{1-n})_s}
\lab{W_hyp} \ee We can present polynomials $\Phi_n(z)$ in the form
\be \Phi_n(z) = \sum_{s=0}^n W_{ns} z^s =\frac{(-1)^n
q^{-n/2}}{1+q^{-n}} \: {_2}{\phi}_1\left( {q^{-n}, -q \atop
-q^{1-n}}; -zq^{-1/2} \right), \lab{P_hyp} \ee where
${_2}{\phi}_1$ stands for the basic hypergeometric function
defined as \cite{KS}
$$
{_2}{\phi}_1\left( {\alpha, \beta \atop \gamma}; z \right) \equiv
\sum_{s=0}^{\infty} \frac{(\alpha)_s (\beta)_s}{(q)_s (\gamma)_s}
z^s
$$
(clearly if $\alpha=q^{-n}$ with positive integer $n$ then we have
a finite sum of $n+1$ terms).

The reflection parameters $a_n$ satisfy the condition $-1 < a_n
<1, \; n=0,1,\dots$, hence there exists a unique weight function
$\rho(\theta)$ on the unit circle such that \be c_n = 1/\cosh(wn)
=\int_{0}^{2\pi} \rho(\theta) \cos(n \theta) d \theta, \quad
n=0,1,2,\dots \lab{hyp_int} \ee
\begin{pr}
The weight function $\rho(\theta)$ for the moments \re{hyp_mom} is
expressed in terms of the Jacobi $\dn$ function: \be \rho(\theta)
= \frac{K(k)}{\pi^2} \: \dn(K(k)\theta/\pi), \lab{hyp_rho} \ee
where $K(k)$ is complete elliptic integral of the first kind and
the elliptic modulus $k$ is determined from the transcendent
equation \be w=\pi K'(k)/K(k) \lab{eq_k} \ee
\end{pr}
Proof of this statement immediately follows from the Fourier
expansion \re{Four_dn} of the $\dn$ function. The function
$\rho(\theta)$ satisfies the required symmetry property
$\rho(2\pi-\theta) = \rho(\theta)$ which follows from the
periodicity property $\dn(u+ 2K)= \dn(u)$. Note that to every
value of $w>0$ there corresponds a unique value of the modulus $k$
in the canonical interval $(0,1)$ such that $k\to 1$ when $w \to
0$ and $k \to 0$ when $w \to \infty$.

The obtained polynomials $P_n(z)$ coincide with a special case of
the family of biorthogonal polynomials introduced by Pastro
\cite{Pas}. Indeed, the Pastro polynomials depend on two
parameters $\alpha, \beta$ and are expressed as
$$
P_n(z) =A_n \:  {_2}{\phi}_1\left( {q^{-n}, q^{\alpha} \atop
q^{2-n-\beta}}; zq^{3/2-\beta} \right)
$$
where $A_n$ is a constant not depending on $z$. Comparing with our
expression \re{P_hyp} we find that $q^{\alpha-1} = q^{\beta-1}=-1$
(we need also to rescale the argument $z \to -z/q$). Thus our
polynomials correspond to a special case of the Pastro polynomials
with fixed values of the parameters $\alpha, \beta$. It is
interesting to note that the weight function $\rho(\theta)$ is
expressed in terms of the elliptic function. This can be also
derived from the results of Pastro \cite{Pas} who obtained
explicit expression of the weight function $\rho(\theta)$ for
arbitrary values of the parameters $\alpha, \beta$. Pastro
presented his result for $\rho(\theta)$ in terms of an infinite
product. In our case this infinite product coincides with
corresponding formula for the function $\dn(x)$.

The "reflected sign" polynomials $\tilde \Phi_n(z) = (-1)^n
\Phi_n(-z)$ have the moments $\tilde c_n = (-1)^n /\cosh(wn)$. It
is easily seen that corresponding weight function on the unit
circle for these polynomials is
$$
{\tilde \rho}(\theta) = \rho(\theta+ \pi) =\frac{k' K(k)}{\pi^2}
\: \frac{1}{\dn(K(k)\theta/\pi)}
$$

\section{Special case of polynomials orthogonal on regular polygons}
\setcounter{equation}{0} So far, we assumed that the parameter $w$
satisfies the restriction $w \ne 4K M_1/M_2$ with some integer
$M_1, M_2$. This restriction is necessary in order for polynomials
$\Phi_n(z)$ to be nondegenerate.

In this section we consider the case when this restriction is
omitted. This means that a degeneration occurs, i.e. $a_{N-1}= \pm
1$ for some positive integer $N$. Assume that under this condition
the polynomial $\Phi_N(z)$ has only simple zeros: \be \Phi_N(z)
=(z-z_0)(z-z_1) \dots (z-z_{N-1}) \lab{Phi_N} \ee i.e. $z_i \ne
z_k$ when $i \ne k$. It can be easily showed \cite{Ger} that under
conditions $|a_n|<1, n=0,1,\dots, N-2$ all these zeros $z_s$ lie
on the unit circle: $|z_s|=1$.

In this case we can consider only a finite set of polynomials
$\Phi_0(z), \Phi_1(z), \dots \Phi_{N-1}(z)$ which are orthogonal
on the finite set of these zeros on the unit circle \cite{Ger} \be
\sum_{s=0}^{N-1} \rho_s \Phi_n(z_s) \Phi_m(1/z_s) = h_n \:
\delta_{nm}, \quad n,m=0,1,\dots, N-1 \lab{finite_ort} \ee where
the discrete weights $\rho_s$ can be expressed as \cite{Zhe2} \be
\rho_s = \frac{h_{N-1}}{\Phi_{N-1}(1/z_s) \Phi_N'(z_s)}
\lab{finite_w} \ee

As an example, consider the case of $cn$ polynomials and condition
\be w=KM/N, \quad M=1,3,\dots, N-2 \lab{cond_N_cn} \ee (it is
assumed that $M$ is co-prime with $N$). From \re{a_cn} we see that
under this condition $a_{2N-1}=-1$. Thus we need to calculate
zeros of the polynomial $\Phi_{2N}(z)$. We cannot substitute
directly the value \re{cond_N_cn} into the explicit expression for
$\Phi_{2N}(z)$ because in this case there is indeterminacy in
expansion coefficients. Nevertheless we can put $w = KM/N +
\epsilon$ and then take the limit $\epsilon =0$ (because from
recurrence relation \re{rec_PC} it is seen that expansion
coefficients of the polynomials $\Phi_{n}(z)$ depend on $w$
continuously). Then we easily obtain that under condition
\re{cond_N_cn} only the last and the first terms survive in the
power expansion of $\Phi_{2N}(z)$:
$$
\Phi_{2N}(z) = z^{2N} +1
$$
and hence all zeros are simple: \be z_s = \exp(\pi i (s-1/2)/N),
\quad s=-N+1, -N+1, \dots N \lab{zer_2N} \ee Note that these zeros
coincide with corresponding spectral points \re{zs_cn} when
$w=KM/N$ and $M$ is co-prime with $N$.

We thus obtained explicitly spectral points $z_s$ in the
finite-dimensional case. In order to get expression for the
weights $\rho_s$ we should calculate the expression \be
\Phi_{2N-1}(1/z_s) = \sum_{m=0}^{2N-1} (-1)^{m+1}
\cn(K(2N-m-1)M/N) E^{2N-1}_m \exp(- i \pi m (s-1/2)/N)
\lab{2N-1_s} \ee However we can avoid direct calculation of the
sum \re{2N-1_s} if we can construct a finite discrete measure on
the unit circle which gives us appropriate values for the moments.
In what follows we will consider only the simplest case when $N$
is prime number and $M=1$. We have the set of $2N$ moments
$$
c_n = \cn(Kn/N), \quad n=0,1,\dots, 2N-1
$$
We should find such values $\rho_j$ that \be c_n = \cn(Kn/N) =
\sum_{j=-N+1}^{N} \rho_j z_j^n = \sum_{j=-N+1}^{N} \rho_j \exp(i
\pi n (j-1/2)/N) \lab{c_n_N} \ee

To do this we return to the Fourier series \re{Four_cn} for the
$\cn$ function and put $z=Kn/N$. It is convenient to present the
summation variable $s$ modulo $2N$, i.e. $s=j+2Nm$, where $j=-N+1,
-N+2, \dots N$ and $m$ can takes all integer values $-\infty < m <
\infty$ . Then formula \re{Four_cn} can be presented in the form
\be \cn(Kn/N) = \frac{\pi}{kK} \sum_{j=-N+1}^{N} e^{i\pi
n(j-1/2)/N} \: S(j;N)  =  \frac{\pi}{kK} \sum_{j=-N+1}^{N} z_j^n
\: S(j;N)   \lab{disc_F} \ee where \be
S(j;N)=\sum_{m=-\infty}^{\infty} \frac{1}{q^{j+2Nm-1/2} +
q^{-j-2Nm+1/2}} \lab{S_def} \ee We thus have \be \rho_j
=\frac{\pi}{kK} S(j;N) \lab{rho_j_N} \ee

In order to calculate the sum $S(j;N)$ explicitly we need the
\begin{lem}
For any $q$ such that $0<q<1$ and any real $\alpha$ introduce the
function  \be F(\alpha;q) = \sum_{s=-\infty}^{\infty}
\frac{1}{q^{n+\alpha} + q^{-n-\alpha}}.  \lab{F_alpha} \ee Then
the function $F(\alpha;q)$ has two equivalent representations:
either \be F(\alpha;q) = \frac{1}{q^{\alpha} + q^{-\alpha}} \:
\frac{(-q^{1+2\alpha},-q^{1-2\alpha},q^2,q^2;q^2)_{\infty}}{(-q^{2+2\alpha},
-q^{2-2\alpha}, q, q;q^2)_{\infty}} \lab{F_prod} \ee or \be
F(\alpha;q) = \frac{K}{\pi} \: \dn(2\alpha K';k'), \lab{F_dn} \ee
where as usual $q=\exp(-\pi K'/K)$. We adopt standard notation for
the q-shifted factorial \cite{GR} $(z;q)_n = (1-z) (1-zq) \dots
(1-zq^{n-1})$ and $(z_1,z_2, \dots, z_N;q)_n$ stands for product
of q-shifted factorials $(z_1;q)_n \dots (z_N;q)_n$.
\end{lem}
The proof of this Lemma follows directly form the Ramanujan
summation formula  for the bilateral ${_1}\psi_1$ basic
hypergeometric function \cite{GR}: \be
 \sum_{s=-\infty}^{\infty} \frac{(c/b;q)_n b^n}{(aq;q)_n} =
 \frac{(c,q/c,abq/c,q;q)_{\infty}}{(b, aq, aq/c,bq/c;q)_{\infty}}
\lab{Ram} \ee where $(x;q)_n$ means q-shifted factorial. In our
case it is sufficient to put $c=ab, \: a=-q^{\alpha}, \: b=
q^{1/2}$. Then the Ramanujan formula \re{Ram} gives
$$
F(\alpha;q^{1/2})= \sum_{s=-\infty}^{\infty}
\frac{1}{q^{(n+\alpha)/2} + q^{-(n+\alpha)/2}} =
\frac{1}{q^{\alpha/2} + q^{-\alpha/2}} \:
\frac{(-q^{1/2+\alpha},-q^{1/2-\alpha},q,q;q)_{\infty}}{(-q^{1+\alpha},
-q^{1-\alpha}, q^{1/2}, q^{1/2};q)_{\infty}}
$$
Then we replace $q \to q^2$ and obtain formula \re{F_prod}.
Formula \re{F_dn} is then obtained  from the well-known
representation of  the Jacobi elliptic functions in terms of
infinite products \cite{WW}. Another method to prove \re{F_dn} is
using the well known Poisson summation formula \cite{Akhiezer2}
\be \sum_{s=-\infty}^{\infty} f(t+sT) = \frac{1}{T}
\sum_{s=-\infty}^{\infty} \hat f\left(\frac{2 \pi s}{T}\right)
\exp(2 i \pi s t/T ) \lab{Poisson} \ee where the functions $f(x),
\hat f(x)$ are connected by the Fourier transform:
$$
\hat f(x) = \int_{-\infty}^{\infty} f(t) e^{-ixt} dt
$$
Applying formula \re{Poisson} to the function $f(x) =
1/\cosh(\beta x)$ with $q=e^{-\beta}$ and using the Fourier
expansion \re{Four_dn} of the dn-function we arrive at formula
\re{F_dn}.

Now we see from \re{S_def} and \re{F_prod} that
$$
S(j;N) = F\left(\frac{j-1/2}{2N}; q^{2N}\right)
$$
and we can express the sum $S(j;N)$ in the one of two equivalent
forms \re{F_prod} or \re{F_dn}. In the "product form" \re{F_prod}
we have \be S(j;N) = \frac{1}{q^{j-1/2} + q^{1/2-j}} \:
\frac{(-q^{2N+2j-1}, -q^{2N-2j+1}, q^{4N}, q^{4N}
  ;q^{4N})_{\infty}}{(-q^{4N+2j-1}, -q^{4N-2j+1}, q^{2N}, q^{2N}  ;q^{4N})_{\infty}}
\lab{prod_S} \ee

In the "$\dn$-form" we have \be S(j;N) = \frac{\tilde K}{\pi} \;
\dn\left(\frac{j+1/2}{N} \tilde K'; \tilde k' \right) \lab{S_dn}
\ee where $\tilde K, \tilde K'$ and $\tilde k, \tilde k'$
correspond to $\tilde q = q^{2N} = \exp(-2\pi N  K'/K)$ instead of
$q$. Such transformation (from $q$ to $q^{2N}$) correspond to
so-called main $2N$-order transformation of Jacobi elliptic
functions \cite{Akhiezer}. In more details, we have
\cite{Akhiezer} \be \tilde K = \frac{K}{2N\mu}, \quad \tilde K' =
\frac{K'}{\mu}, \quad \tilde k = k^{2N} \: \prod_{r=1}^N
\sn^4\left(\frac{(2r-1)K}{2N};k \right) \lab{til_Kk} \ee where
$$
\mu = \prod_{r=1}^N \frac{\sn^2\left(\frac{(2r-1)K}{2N};k
\right)}{\sn^2\left(\frac{rK}{N};k \right)}
$$
It is interesting to note that the $\dn$-function in \re{S_dn} can
be explicitly expressed as a rational function
$$
S(j;N) = R_N(\phi(j))
$$
of the function
$$
\phi(j) = \dn^2\left(\frac{(j+1/2)K'}{N};k'  \right)
$$
We will not consider this form here (see, e.g. \cite{Akhiezer} for
details, where relation with so-called Zolotarev rational
functions of the best approximation is discussed).

We thus have the following
\begin{pr}
Assume that $w=K/N$, where $N=1, 2,3,\dots$ . Define the
reflection parameters
$$
a_{2n} = \cn(w(2n+1)), \quad \quad a_{2n+1}= - \dn(w(2n+2)), \quad
n=0,1,\dots, N-1
$$
such that $a_{2N-1}=-1$. Define the monic polynomials $\Phi_n(z),
n=0,1,\dots 2N$ by recurrence relation \re{Sz_real}.

Then the polynomial $\Phi_{2N}(z)=z^{2N}+1$ has simple zeros
$z_s=\exp(\pi i (s-1/2)/N), \quad s=-N+1, -N+1, \dots N$. The
polynomials $\Phi_n(z), \; n=0,1,\dots, 2N-1$
 are orthogonal on the finite set
$z_s$ of spectral points on the regular $2N$-gon: \be
\sum_{j=-N+1}^N \rho_j \Phi_n(z_j) \Phi_m(1/z_j) = h_n \:
\delta_{nm} \lab{ort_fin_cn} \ee where the discrete weights
$\rho_j$ are given by formula \re{rho_j_N} with $S(j;N)$ given by
\re{prod_S} or \re{S_dn}.
\end{pr}
{\it Remark}. One can consider the generic case $w=KM/N$ with
arbitrary co-prime integers $M,N$. This will lead to slightly
modified formulas for the weights $\rho_s$. The case of
$\dn$-circle polynomials can be considered in a similar way.

\section{Polynomials orthogonal on the interval}
\setcounter{equation}{0} Now we can relate with our polynomials
orthogonal on the unit circle $P_n(z)$ some polynomials $S_n(x)$
orthogonal on an interval of the real axis. This relation was
first proposed by Delsarte and Genin \cite{DG1} (see also
\cite{Zhe1} for further results concerning the Delsarte-Genin
transformation (DGT)).

Recall the main properties of the DGT. Given polynomials $P_n(z)$
with real reflection parameters $a_n$ satisfying the restriction
$-1 < a_n <1$ we construct the polynomials $S_n(x)$ by the formula
\be S_n(x) = \frac{z^{-n/2}\Phi_n(z) +  z^{n/2}\Phi_n(1/z)
}{2^n(1-a_{n-1})}, \lab{S_P} \ee where $x=(z^{1/2} + z^{-1/2})/2$
and it is assumed that $z^{\pm 1/2} = r^{1/2} e^{\pm i \theta/2}$
when $z= r e^{i \theta}$. It is easily verified that the
polynomials $S_n(x)$ are monic $n$-th degree polynomials in $x$:
$S_n(x) = x^n + O(x^{n-1})$. Moreover, the polynomials $S_n(x)$
are {\it symmetric}, i.e. $S_n(-x) = (-1)^n S_n(x)$ and satisfy
the recurrence relation \be S_{n+1}(x) + v_n S_{n-1}(x) = x S_n(x)
\lab{rec_S} \ee where \be v_n = \frac{1}{4}(1+a_{n-1})(1-a_{n-2})
\lab{u_a} \ee Note that for $n=0$ we have formally $v_0=0$ due to
assumed condition $a_{-1}=-1$. Due to the same condition we have
also $v_1 = (1+a_0)/2$.

Under restriction $-1 < a_n <1$ the recurrence coefficients $v_n$
are strictly positive $v_n>0, n=1,2,\dots$. Thus the polynomials
$S_n(x)$ are orthogonal on the interval $[-1,1]$ with respect to a
positive weight function \be \int_{-1}^1 S_n(x) S_m(x) w(x) dx =
\kappa_n \: \delta_{nm} \lab{ort_S} \ee where \be w(x) =
\frac{\rho(\theta)}{\sin(\theta/2)}, \quad x=\cos(\theta/2)
\lab{w_S} \ee and
$$
\kappa_n = v_1 v_2 \dots v_n = 2^{1-2n} h_n(1-a_{n-1}^2)^{-1}
$$
The orthogonality relation can be presented in an equivalent form
\be \int_{-1}^1 w(x) \frac{S_n(x) S_m(x) + S_n(-x) S_m(-x)}{2} dx
= \kappa_n \: \delta_{nm} \lab{ort_S1} \ee which follows from the
property $w(-x) = w(x)$ for all symmetric polynomials $S_n(x)$.
Relation \re{ort_S1} is more convenient when the measure is
discrete one.

Introduce the moments $$M_n= \int_{-1}^1 w(x) x^n d x, \quad
n=0,1,2,\dots$$ corresponding to the weight function $w(x)$ of the
polynomials $S_n(x)$. It is easily seen  that the moments $M_n$
are expressed in terms of the moments $c_n$ as \cite{DG1} \be M_n
= 2^{-n} \: \sum_{j=0}^n {n \choose j} c_{j-n/2} \lab{M_c} \ee
when $n$ is even and $M_n =0$ when $n$ is odd.

In the case of $cn$-polynomials we obtain that corresponding
symmetric polynomials $S_n(x)$ are orthogonal \be
\sum_{s=-\infty}^{\infty} \frac{S_n(x_s)S_m(x_s) +
S_n(-x_s)S_m(-x_s)}{2} \rho_s = \kappa_n \:\delta_{nm}
\lab{cn_S_ort} \ee on the dense set of the points
$$
x_s=\cos\left(\frac{\pi w(s-1/2)}{2K}\right)
$$
with the weights $\rho_s$ is given by \re{W_cn}.

In the case of $dn$-polynomials we obtain similar orthogonality
relation \re{cn_S_ort} but now the orthogonality grid is
$$
x_s=\cos\left(\frac{\pi ws}{2K}\right)
$$
and the discrete weights $\rho_s$ are given by \re{W_dn}.

As expected, both $cn$ and $dn$-type polynomials $S_n(x)$ are
orthogonal with respect to measures which are positive and dense
on the interval $[-1,1]$.

It is well known \cite{Chi} that to any family of symmetric
orthogonal polynomials $S_n(x)$ one can associate two families of
orthogonal polynomials $P_n(x)$ and $Q_n(x)$ by the formulas
$S_{2n}(x)= 2^{-n}P_n(2x^2-1)$ and $S_{2n+1}(x) =2^{-n}x
Q_n(2x^2-1) $ . Polynomials $P_n(y)$ and $Q_n(y)$ are monic
$P_n(y) = y^n + O(y^{n-1}), \; Q_n(y) = y^n + O(y^{n-1})$ and
orthogonal on the same interval $[-1,1]$: \be \int_{-1}^1 P_n(y)
P_m(y) w(\sqrt{(y+1)/2})(y+1)^{-1/2} dy =2^{2n-1} \kappa_{2n} \:
\delta_{nm} \lab{ort_Py} \ee  and  \be \int_{-1}^1 Q_n(y) Q_m(y)
w(\sqrt{(y+1)/2})(y+1)^{1/2} dy =2^{2n-1} \kappa_{2n+1} \:
\delta_{nm} \lab{ort_Qy} \ee

The polynomials $P_n(x)$ satisfy the recurrence relation
\cite{Chi} \be P_{n+1}(x) + b_n P_n(x) + u_n P_{n-1}(x) = xP_n(x),
\lab{rec_P} \ee where \be u_n = 4v_{2n} v_{2n-1}, \quad b_n =
2(v_{2n} + v_{2n+1})-1 \lab{ub_vP} \ee Analogously, polynomials
$Q_n(x)$ satisfy the recurrence relation \be Q_{n+1}(x) + b_n
Q_n(x) + u_n Q_{n-1}(x) = xQ_n(x), \lab{rec_Q} \ee where \be u_n =
4v_{2n} v_{2n+1}, \quad b_n = 2(v_{2n+2} + v_{2n+1})-1 \lab{ub_vQ}
\ee

Using the representation \re{u_a} we obtain explicit formulas for
the recurrence coefficients of the polynomials $P_n(x)$:  \be u_n
= \frac{1}{4} (1+a_{2n-1})(1-a^2_{2n-2})(1-a_{2n-3}), \; b_n =
\frac{1}{2}\left( a_{2n}(1-a_{2n-1}) -a_{2n-2}(1+a_{2n-1}) \right)
\lab{ub_aP} \ee Similarly, for the recurrence coefficients of the
polynomials $Q_n(x)$ we have \be u_n = \frac{1}{4}
(1+a_{2n})(1-a^2_{2n-1})(1-a_{2n-2}), \; b_n = \frac{1}{2}\left(
a_{2n+1}(1-a_{2n}) -a_{2n-1}(1+a_{2n}) \right) \lab{ub_aQ} \ee

 Formulas \re{ub_aP}, \re{ub_aQ} were first
obtained by Geronimus \cite{Ger} (see also \cite{Simon}).

Explicitly polynomials $P_n(x)$ are expressed in terms of
$\Phi_n(z)$ as \be P_n(x) = \frac{z^{-n} \Phi_{2n}(z) + z^n
\Phi_{2n}(1/z)}{2^{n} (1-a_{2n-1})}, \quad x=(z+z^{-1})/2
\lab{Sz_P} \ee (formula \re{Sz_P} is due to Szeg\H{o}). Similarly
for polynomials $Q_n(x)$ we have \be Q_n(x) = \frac{z^{-n-1/2}
\Phi_{2n+1}(z) + z^{n+1/2} \Phi_{2n+1}(1/z)}{2^{n} (z^{1/2} +
z^{-1/2}) (1-a_{2n})}, \quad x=(z+z^{-1})/2 \lab{Sz_Q} \ee

In what follows we will consider only the polynomials $P_n(x)$.
All formulas for the companion polynomials $Q_n(x)$ are obtained
in a similar way.

Polynomials $P_n(x)$ are orthogonal on the interval $[-1,1]$ with
respect to the weight function \be w(x) =
\frac{\rho(\theta)}{|\sin \theta|} =
\frac{\rho(\theta)}{\sqrt{1-x^2}} \lab{weight_P} \ee

Assume that polynomials $\Phi_n(z)$ have the expression
$$
\Phi_n(z) = \sum_{s=0}^{n} W_{ns} z^s
$$
with some coefficients $W_{ns}$. Then from \re{Sz_P} we obtain
corresponding expansion formula for polynomials $P_n(x)$
\cite{Atk}: \be P_n(x) = 2^{1-n} (1-a_{2n-1})^{-1} \left(W_{2n,n}
+ \sum_{s=1}^n (W_{2n,n+s} + W_{2n,n-s})T_s(x) \right),
\lab{P_Cheb} \ee where $T_n(x) = \cos(\theta n)$
 is the Chebyshev polynomial of the first kind (it is assumed that
 $x=\cos\theta$).

Thus if polynomials $\Phi_n(z)$ on the unit circle have explicit
expression as expansion in terms of monomials $z^n$ then the
corresponding polynomials $P_n(x)$ on the interval have explicit
expression as expansion in terms of the Chebyshev polynomials
$T_n(x)$.

Return to our elliptic $\cn$-polynomials. Corresponding
polynomials $P_n(x)$ on the interval have the explicit recurrence
coefficients $$u_n= \frac{\sn^2(w(2n-1)) (1-\dn(2wn)
(1+\dn(2w(n-1)) }{4}$$
$$
b_n = \frac{\cn(w(2n+1))(1 + \dn(2wn)) - \cn(w(2n-1))(1 -
\dn(2wn))}{2}
$$

%To simplify expressions for the recurrence coefficients we will
%use well known formulas \cite{WW}
%$$
%\frac{1-\cn u}{\sn u} = \frac{\sn(u/2) \dn(u/2)}{\cn(u/2)}, \quad
%\frac{1+\cn u}{\sn u} = \frac{\cn(u/2)}{\sn(u/2) \dn(u/2)}, \quad
%\frac{1-\dn u}{\sn u} = k^2\: \frac{\sn(u/2) \cn(u/2)}{\dn(u/2)},
%\quad \frac{1+\dn u}{\sn u} = \frac{\cn(u/2) \dn(u/2)}{\sn(u/2)}
%$$

Using explicit expression \re{exp_cn_even} for the expansion
coefficients $W_{ns}$ we obtain in case of the $cn$ polynomials
\be P_n(x) = (-1)^n 2^{1-n} (1+ \dn(2wn))^{-1} \left( \dn(wn)
E^{2n}_n + \sum_{s=1}^n (-1)^{s} E^{2n}_{n+s}(\dn(w(n-s)) +
\dn(w(n+s))T_s(x) \right) \lab{P_cn} \ee These polynomials satisfy
the orthogonality relation \be \sum_{s=-\infty}^{\infty} \rho_s
P_n(x_s)P_m(x_s) = H_n \delta_{nm}, \quad H_n = u_1 u_2 \dots u_n
\lab{ort_P_cn} \ee where the spectral points are
$$
x_s= \cos\left(\frac{\pi w (s-1/2)}{K} \right)
$$
and the weights $\rho_s$ are given by \re{W_cn}.

The polynomials $Q_n(x)$ defined by \re{Sz_Q} are orthogonal on
the same set of spectral points but the measure acquires an
additional linear factor $1+x$:
$$
\sum_{s=-\infty}^{\infty} \rho_s (1+x_s) Q_n(x_s)Q_m(x_s) =0,
\quad \mbox{if} \quad n \ne m
$$

Clearly, we have the case of a dense point spetrum on the
interval $[-1,1]$.

In the limit $k=1$ we obtain the recurrence coefficients for the
polynomials $P_n(x)$
$$
u_n=\frac{1}{4} \: \frac{(1-q^n)^2(1-q^{2n-1})^2 (1+q^{n-1})^2
}{(1+q^{2n-1})^2 (1+q^{2n}) (1+q^{2n-2})}
$$
$$
b_n = q^{n-1/2} \: \frac{2(q+1) q^n -(1-q)(1-q^{2n})
}{(1+q^{2n-1})(1+q^{2n+1})}
$$
where $q=e^{-2w}$ (of course, this $q$ should not be confused with
already introduced "elliptic" $q$ by \re{qKK}).

It is easily verified that these recurrence coefficients
correspond to the Askey-Wilson orthogonal polynomials \cite{AW},
\cite{KS} $p_n(z;a,b,c,d|q)$ with the parameters $a=c=q^{1/2},
b=-d=1$. When $max(|a|,|b|,|c|,|d|)\le 1$ then the Askey-Wilson
polynomials are orthogonal on the interval $[-1,1]$ with respect
to a positive continuous measure \cite{AW}, \cite{KS}
$$
\int_{-1}^1 \frac{w(x;a,b,c,d|q)}{\sqrt{1-x^2}} p_n(z;a,b,c,d|q)
p_m(z;a,b,c,d|q) dx = h_n \: \delta_{nm}
$$
where
$$
w(x;a,b,c,d|q)= \frac{h(x;1)
h(x;-1)h(x;q^{1/2})h(x;-q^{1/2})}{h(x;a)h(x;b)h(x;c)h(x;d)}
$$
and
$$
h(x;\alpha) = \prod_{s=0}^{\infty} [1-2\alpha x q^s + \alpha^2
q^{2s}]
$$
Thus in our case the weight function becomes \be w(x) =
\frac{h(x;-q^{1/2})}{h(x;q^{1/2})} = \prod_{s=0}^{\infty} \frac{1+
2 x q^{s+1/2} + q^{2s+1}}{1- 2 x q^{s+1/2} + q^{2s+1}}
\lab{w_AW_d} \ee Expression \re{w_AW_d} can be compared with the
expression of the Jacobi dn-elliptic function in terms of the
infinite product \cite{WW} \be \dn(u;k) =\sqrt{k'} \:
\prod_{s=0}^{\infty} \frac{1+ 2 \cos(2\theta) q^{2s+1} +
q^{4s+2}}{1- 2 \cos(2\theta) q^{2s+1} + q^{4s+2}}, \quad u = 2K
\frac{\theta}{\pi} \lab{dn_prod} \ee We see that the weight
function $w(x)$ for the corresponding Askey-Wilson polynomials is
proportional to the elliptic $dn(u)$ function under the
identification $q \to q^2, \: x \to \cos(2\theta)$.

For the case of $Q$-polynomials defined by recurrence coefficients
\re{ub_aQ} we have in the limit $k =1$ again the Askey-Wilson
polynomials $p_n(x;a,b,c,d|q)$ with the parameters $a=c=q^{1/2},
b=1, d=-q$. For the weight function we obtain
$$
w(x) = \frac{h(x;-1)h(x;-q^{1/2})}{h(x;-q) h(x;q^{1/2})} = 2(1+x)
\: \frac{h(x;-q^{1/2})}{h(x;q^{1/2})}
$$
which differs from the weight function \re{w_AW_d} only by a
linear factor $2(1+x)$.

It is interesting to note that despite the hyperbolic limit $k \to
1$ some "smell" of the elliptic functions remains in the
expression for the weight function of the Askey-Wilson
polynomials. Askey and Wilson themselves discussed intriguing
cases when the weight function $w(x)$ is expressible in terms of
the Jacobi elliptic functions $\sn(u), \cn(u), \dn(u)$ \cite{AW}.
We see that the case of $dn$-weight function of the Askey-Wilson
polynomials has a natural "elliptic" generalization leading to
orthogonal polynomials with dense measure on the same interval
$[-1,1]$. It would be interesting to search other possible
elliptic generalizations of the Askey-Wilson polynomials.

\section{A general scheme leading to a dense point spetrum on the unit circle}
\setcounter{equation}{0} The examples of polynomials on the unit
circle constructed in the previous sections can be generalized. In
this section we propose a wide class of POC having a positive
dense point spetrum on the unit circle.

Let $f(x)$ be a continuous function of a real variable $x$ with
three main properties:

(i) $f(x)$ is an even function $f(-x)=f(x)$ normalized by the
condition $f(0)=1$;

(ii) function $f(x)$ is periodic with some real period $T$:
$f(x+T) = f(x)$;

(iii) the function $f(x)$ has the Fourier expansion \be f(x) =
\sum_{n=-\infty}^{\infty} A_n \: e^{2\pi i n x/T}, \quad
A_{-n}=A_n, \lab{Four_f} \ee such that all the Fourier
coefficients are nonnegative: $A_n \ge 0, \; n=0,1,2,\dots$.
Moreover, we assume that there are infinity many positive
coefficients $A_n >0$ (otherwise the problem is trivial). Of
course, we assume that the Fourier series converges everywhere on
the real axis.

We define the moments by the formula \be c_n = f(wn), \quad
n=0,1,2,\dots, \lab{mom_f} \ee where $w$ is an arbitrary positive
parameter such that \be w \ne TM/N \lab{non_NM} \ee for some
integers $M,N$. Then it is obvious that $c_{-n} = c_n$ and the
moments $c_n, \; n=0,1,2,\dots$ are all distinct: $c_n \ne c_m$ if
$n \ne m$. Using these moments $c_n$ we can construct the Toeplitz
determinants $\Delta_n$ and corresponding OPC by \re{deterP}.

\begin{th}
Under conditions (i)-(iii) the Toeplitz  determinants $\Delta_n$
are all positive $\Delta_n>0$, or equivalently all reflection
parameters satisfy the restriction $-1 < a_n <1 $. The polynomials
$P_n(z)$ are orthogonal on the unit circle with respect to a
positive dense point spetrum \be \sum_{s=-\infty}^{\infty}
\rho_s \Phi_n(z_s)\Phi_m(1/z_s) = h_n \: \delta_{nm}
\lab{sing_ort} \ee where the orthogonality grid on the unit circle
is \be z_s = e^{2 \pi i w s/T} \lab{grid_zs} \ee and the discrete
weights coincide with the Fourier coefficients \be \rho_s = A_s,
\; s=\pm 1, \pm 2, \dots \lab{rho_s} \ee
\end{th}
{\it Proof}. From the Fourier expansion \re{Four_f} we have the
representation for the moments in the form
$$
c_n = f(wn) = \sum_{s=-\infty}^{\infty} A_s z_s^n =
\sum_{s=-\infty}^{\infty} \rho_s z_s^n = \int_C z^n d \mu(z), \;
n=0, \pm 1, \pm 2, \dots
$$
where $z_s$ is given by \re{grid_zs} and $d \mu(z)$ is the measure
on the unit circle corresponding to the positive discrete weights
$\rho_s>0$. Whence the polynomials $P_n(z)$ possess orthogonality
property \re{sing_ort}. The measure $\mu(z)$ is positive and
dense on the unit circle. Indeed, positivity follows
from the property $A_n > 0$ for infinitely many values of $n$. On
the other hand we see that in any small vicinity of the point
$z_s$ on the unit circle there exists at least one point $z_{s'}$
from the set of spectral points $z_s$ due to condition
\re{non_NM}. Thus the points $z_s$ form a Cantor dense set on the
unit circle. This means that the corresponding point spectrum is dense.

Now from general properties of the polynomials orthogonal on the
unit circle we conclude that all $\Delta_n$ are positive
$\Delta_n>0$ and equivalently, that the reflection parameters
satisfy the restriction $-1 < a_n < 1$.

We thus see that starting from arbitrary even periodic continuous
function $f(x)$ having infinitely many positive Fourier
coefficients $A_n$ we can construct corresponding OPC with
a positive dense point spetrum on the unit circle.

The Toeplitz determinants have the expression \be \Delta_n = \left
|
\begin{array}{cccc} 1 & f(w) & \dots &
f(w(n-1))\\ f(w) & 1 & \dots & f(w(n-2)) \\ \dots & \dots & \dots & \dots\\
f(w(n-1)) & f(w(n-2)) & \dots & 1 \end{array} \right |
\lab{Toepl_w} \ee

All the reflection coefficients of these polynomials will satisfy
the restriction $-1 < a_n <1$. Of course, if the function $f(x)$
and its Fourier coefficients $A_n$ are known explicitly, we have
explicit orthogonality relation \re{sing_ort}. However, explicit
expression for the reflection coefficients $a_n$ can be obtained
only for some exceptional cases.

As far as we know the first (and so far the only) explicit example
of such type was proposed by Magnus \cite{Mag}, \cite{Mag1} who
considered the function $f(x)$ defined as
$$
f(x) = (-1)^{[x] +1}(2(x -  [x]) +1),
$$
where $[x]$ means integer part, i.e. the greatest integer less
than or equal to x. Equivalently, the function $f(x)$ can be
defined as $f(x)=1-2|x|$ for $-1 \le x \le 1$ and then continued
for the whole real axis by 2-periodicity.

Clearly the function $f(x)$ is even $f(-x)=f(x)$ and  periodic
with period 2: $f(x+2)=f(x)$. Moreover, there is an obvious
property $f(x+1) =- f(x)$. If two numbers $x$ and $y$ have the
same integer part $[x]=[y]$ then obviously \be \frac{f(x) -
f(y)}{x-y} = -(-1)^{[x]} \lab{int_f} \ee

The function $f(x)$ has the Fourier expansion
$$
f(x) = \frac{4}{\pi^2} \: \sum_{s=-\infty}^{\infty} \frac{e^{\pi i
 x (2s+1)}}{(2s+1)^2}
$$
Thus $A_n = 4/(\pi^2 n^2)$ if $n$ is odd and $A_n=0$ if $n$ is
even. We see that all $A_n$ are nonnegative, hence conditions for
existence of a positive dense point spetrum are fulfilled
and corresponding polynomials are orthogonal on the unit circle
\be \sum_{s=-\infty}^{\infty} (2s+1)^{-2} \Phi_n(z_s)\Phi_m(1/z_s)
= h_n \: \delta_{nm} \lab{sing_ort_Mag} \ee with some positive
normalization coefficients $h_n$. The spectral points on the unit
circle are
$$
z_s = e^{i \pi(2s+1)w}
$$
Magnus showed \cite{Mag1} that the polynomials $\Phi_n(z)$ (and
corresponding reflected parameters $a_n$) can be found explicitly.
This explicit expression is closely related with the expansion of
the irrational number $w$ to the continued fraction \be w = q_0 +
{1\over\displaystyle q_1 + {\strut 1 \over\displaystyle q_2 +
{\strut 1 \over {q_3 + \dots }}}}, \qquad \lab{cont_w} \ee
Introduce corresponding convergents $P_i/Q_i, \: i=0,1,2,\dots$,
where $P_0=[w], Q_0=1$ is "zero-step" rational (integer in this
case) approximation of the number $w$ and the rational number
$P_n/Q_n$ is obtained when truncating the continued fraction
\re{cont_w} at the $n$-th step:
$$
\frac{P_n}{Q_n} = q_0 + {1\over\displaystyle q_1 + {\strut 1
\over\displaystyle q_2 + \dots + {\strut 1 \over {q_n}}}}
$$
The $n$-th convergent $P_n/Q_n$ provides the best rational
approximation of the irrational number in the following sense:
$$
|Q_n w - P_n| < |y w- x|
$$
for all possible integers $x,y$ with the restriction $y<Q_n$.

For any given positive integer $n$ we define a sequence of pairs
of nonnegative integers $(n_i, m_i)$ in the following way. Let the
pair $(n_i, m_i)$ provides the $i$-th order best rational
approximation of the number $w$ under condition $1\le n_i \le n$.
In more details this means the following. Define positive
(irrational) numbers $y_i = |w n_i - m_i|$. The pairs $(n_i, m_i),
\; i=1,2,\dots, n$ correspond to the ordering $0<y_1<y_2<y_3 <
\dots < y_n$, where $1 \le n_i \le n$, i.e. the pair $(n_1, m_1)$
yields the best rational approximation to the number $w$ (with
denominators not exceeding $n$), the pair $(n_2, m_2)$ yields then
second (next) best approximation etc. It is well known that the
best approximation $(n_1, m_1)$ is given by an appropriate
continued fraction of the number $w$: $m_1=P_k, \; n_1 =Q_k$,
where $Q_k \le n$. The second approximation $(n_2, m_2)$ as well
all other approximations $(n_i, m_i)$ can be easily expressed in
terms of so-called intermediate fractions (see e.g.
\cite{Khinchin}). Now we can always present the polynomial
$\Phi_n(z)$ in the following form \be \Phi_n(z) = z^n +
G_n^{(1)}z^{n-n_1} + G_n^{(2)}z^{n-n_2} + \dots G_n^{(s)}z^{n-n_s}
+ \dots + G_n^{(n)} \lab{exp_P_A} \ee with some coefficients
$G_n^{(s)}$, where $n_1, n_2, \dots$ correspond to the succeeding
best approximations of the number $w$. (Clearly, the set $n_1,
n_2, \dots n_n$ is a permutation of the set $1,2,\dots ,n$, i.e.
$n_i \ne n_k$ if $i \ne k$ due to irrationality of $w$). It can be
easily showed (see \cite{Mag1} for details) that only 3 first
terms survive in \re{exp_P_A}, i.e. we have \be \Phi_n(z) = z^n +
G_n^{(1)}z^{n-n_1} + G_n^{(2)}z^{n-n_2} \lab{P_3} \ee where the
coefficients $G_n^{(1)}, G_n^{(2)}$ can be found explicitly if
$(n_1, m_1)$ and $(n_2, m_2)$ are known \cite{Mag1}.

Note that in general we have $a_n=0$ for "almost all" $n$. The
nonzero coefficients $a_n$ appear only if for some $n$ we have
$n_1=n$ or $n_2=n$. We thus need information about the best
approximations of the number $w$. Unfortunately for almost all
irrational numbers $w$ it is impossible to present some "explicit"
formula for $n_1, n_2$ (expressing $n_1, n_2$ as an analytic
function of $n$).

In contrast, in our examples of $cn$- and $dn$- elliptic
polynomials we have purely explicit presentation for polynomials
$\Phi_n(z)$, moments $c_n$ and reflection parameters $a_n$.

Is it possible to construct other explicit examples of the
polynomials $\Phi_n(z)$ corresponding to functions $f(x)$
satisfying properties (i)-(iii)?

\section{Remarks concerning "elliptic hypergeometric functions"}
\setcounter{equation}{0} Consider $cn$  or $dn$ polynomials
$\Phi_n(z)$ on the unit circle as some "elliptic hypergeometric
functions" depending on the argument $z$ and discrete
parameter $n$. The parameter $w$ plays the role of the
"deformation" parameter: when $w =0$ we have from \re{P_W} \be
\Phi_n(z) = (z-1)^n \lab{w=0_Phi} \ee for both $cn$ and $dn$
polynomials. Thus in the limit $w=0$ we have the simplest
hypergeometric function
$$
\lim_{w\to 0} \Phi_n(z)= (-1)^n \: {_1}F_0(-n;z) = (z-1)^n.
$$
This allows us to consider the functions $\Phi_n^{(C)}(z)$ and $\Phi_n^{(D)}(z)$ as 
"elliptic deformations" of the hypergeometric function
${_1}F_0(-n;z)$.

Note that from the explicit expression \re{P_W} it is sen that
functions $\Phi_n^{(C)}(z)$ and $\Phi_n^{(D)}(z)$ {\it do not} belong to a class of elliptic
hypergeometric functions ${_{12}}V_{11}$ introduced by Frenkel and Turaev
\cite{FT} and studied intensively during last years (see e.g.
\cite{S_Bailey}, \cite{GR} and many others). In
particular, these functions perhaps do not possess
simple transformation formulas like Bailey transform
\cite{S_Bailey}. Nevertheless, functions $\Phi_n^{(C)}(z)$ and $\Phi_n^{(D)}(z)$  have several
nice properties which are very close to the classical ones. In
particular the recurrence relation \re{3UC} can be considered as a
contiguous relation for 3 elliptic hypergeometric functions
$\Phi_n(z)$ with parameters $n, \: n\pm 1$. The transformation
properties \re{EPP} with respect to the operator ${\cal E}$ can be
considered as an elliptic analogue of the derivatives of the
hypergeometric functions. Finally, the functions $\Phi_n^{(C)}(z)$ and $\Phi_n^{(D)}(z)$ 
possess nice summation formulas for two values of the argument
$z=\pm 1$. Indeed, from \re{Phi1} we have, e.g. for the
$cn$-polynomials
$$
\Phi_{2n}^{(C)}(1) = (1-\cn(w))(1+\dn(2w))\dots (1-\cn(w(2n-1))(1+
\dn(2wn)
$$
and
$$
\Phi_{2n+1}^{(C)}(1) = (1-\cn(w))(1+\dn(2w))\dots (1+
\dn(2wn)(1-\cn(w(2n+1))
$$
(similar expressions hold for $\Phi_{2n}^{(D)}(\pm 1)$).

Thus there exists a much simpler class of "elliptic hypergeometric
functions" with nice properties similar to the classical ones.

{\bf {\Large Acknowledgments}}

The author thanks B.Simon, L.Golinskii, A.Magnus  and S.Tsujimoto for discussions of 
results of the papers.

\bb{99}

\bi{Akhiezer} N.I. Akhiezer, {\it Elements of the Theory of
Elliptic Functions}, 2nd edition, "Nauka", Moscow, 1970.
Translations Math. Monographs {\bf 79}, AMS, Providence, 1990.

\bi{Akhiezer2} N.I.Achieser, {\it Theory of Approximation}, Dover,
1992.

\bi{And} G.E.Andrews, {The Theory of Partitions}, Addison-Wesley,
1976.

%\bi{biask} R. Askey, Comments in: "Gabor Szeg\H{o}:Collected
%Papers", Birkh\"auser, Basel, 1982, v.1, 303-305.

\bi{AW} R.~Askey and J.~Wilson, {\it Some basic hypergeometric
orthogonal polynomials that generalize Jacobi polynomials}, Mem.
Amer. Math. Soc. {\bf 54}, No. 319, (1985), 1-55.

\bi{Atk} F.V.Atkinson, {\it Discrete and Continuous Eigenvalue
Problems}, Academic Press, NY, London, 1964.

\bi{Chi} T. Chihara, {\it An Introduction to Orthogonal
Polynomials}, Gordon and Breach, NY, 1978.

\bi{DG1}  P. Delsarte and Y. Genin, {\it The split Levinson
algorithm}, IEEE Trans. Acoust. Speech Signal Process. {\bf 34}
(1986), 470--478.

\bi{DG2} P. Delsarte and Y. Genin, {\it Tridiagonal approach to
the algebraic environment of Toeplitz matrices I. Basic results},
SIAM J. Matrix Anal. Appl. {\bf 12}, (1991), 220--238.

\bi{DG3} P. Delsarte and Y. Genin, {\it Tridiagonal approach to
the algebraic environment of Toeplitz matrices. II. Zero and
eigenvalue problems}, SIAM J. Matrix Anal. Appl. {\bf 12}, No. 3
(1991), 432--448.

\bibitem{FT} I.B. Frenkel and V.G. Turaev, {\it Elliptic solutions of the
Yang-Baxter equation and modular hypergeometric functions}, The
Arnold-Gelfand Mathematical Seminars, Birkh\"auser, 1997, pp.
171-204.

\bi{GR} G. Gasper and M. Rahman, Basic Hypergeometric Series, 2-nd
Edition, Cambridge University Press, Cambridge, 2006.

\bi{Ger} Ya.L. Geronimus,\quad {\it Polynomials Orthogonal on a
Circle and their Applications}, \\ Am.Math.Transl.,Ser.1 {\bf
3}(1962), 1-78.

%\bi{Ger1} Ya.L.Geronimus, {\it On polynomials orthogonal with
%respect to the given numerical sequence and on Hahn's theorem},
%Izv.Akad.Nauk, {\bf 4} (1940), 215-228 (in Russian).

%\bi{Hahn} W.Hahn, {\it \"Uber die Jacobischen Polynome und Zwei verwandte
%Polynomklassen}, Math.Z. {\bf 39} (1935), 634-638.

%\bi{HN} E.Hendriksen and O. Njastad, {\it Biorthogonal Laurent polynomials
%with biorthogonal derivatives}, Rocky Mount. J. Math. {\bf 21} (1991),
%391-317.

\bi{HR} E.Hendriksen and H. van Rossum, {\it Orthogonal Laurent polynomials},
Indag. Math. (Ser. A) {\bf 89}, 17-36.

\bi{IsMas} M.E.H. Ismail and D. Masson, {\it Generalized orthogonality and
continued fractions}, \\J.Approx.Theory {\bf 83} (1996), 1-40.

%\bi{JT} W.B.Jones and W.J.Thron, {\it Survey of continued fraction methods of
%solving moment problems} in: analytic Theory of Continued Fractions, LNM 932,
%Springer, Berlin, Heidelberg, New York (1981).

\bi{Khinchin} A.Ya.Khinchin, {\it Continued fractions}, Dover,
1997.

\bi{KS} R. Koekoek  and R.F. Swarttouw, {\it The Askey scheme of
hypergeometric orthogonal polynomials and its q-analogue}, Report
94-05, Faculty of Technical Mathematics and Informatics, Delft
University of technology, 1994.

\bi{Mag} A.P. Magnus, {\it Special nonuniform lattice (snul)
orthogonal polynomials on discrete dense sets of points}, J. Comp.
Appl. Math. {\bf 65} (1995), 253--265.

\bi{Mag1} A.P.Magnus, {\it Semi-classical orthogonal polynomials
on the unit circle}. MAPA 3072A Special topics in approximation
theory   http://www.math.ucl.ac.be/membres/magnus/num3/m3xxx99.pdf

%\bi{Mar} P.Maroni, {\it Variations around classical orthogonal
%polynomials. Connected problems}, J.Comp.Appl.Math. {\bf 48}
%(1993), 133-155.

%\bi{Mass} D.Masson, {\it Difference equations, continued fractions, Jacobi
%matrices and orthogonal polynomials}, in Nonlinear Numerical methods and
%Rational Approximations, (A.Cuyt, ed.), 239-257. Reidel Publ.Co. 1988.

\bi{Pas} P.I. Pastro, {\it Orthogonal polynomials and some q-beta integrals
of Ramanujan},\\  J.Math.Anal.Appl. {\bf 112} (1985), 517-540.

\bi{Simon} B.Simon, {\it Orthogonal Polynomials On The Unit
Circle}, AMS, 2005.

\bi{S_Bailey} V.Spiridonov, {\it Elliptic hypergeometric functions}, a complement to the book by G.E. Andrews, R. Askey, and R. Roy, Special Functions, Encycl. of Math. Appl. {\bf 71}, Cambridge Univ. Press, 1999, 
written for its Russian edition. http://xxx.lanl.gov/abs/0704.3099

\bibitem{SZ1} V.Spiridonov and A.Zhedanov, {\it Spectral transformation chains
and some new biorthogonal rational functions}, Commun. Math. Phys.
{\bf 210} (2000), 49-83.

\bi{Sz} G. Szeg\H{o}, Orthogonal Polynomials, fourth edition, AMS,
1975.

\bi{WW} E.T. Whittacker, G.N. Watson, {\em A Course of Modern
Analysis}, Cambridge, 1927.

\bi{Zhe1} A.Zhedanov, {\it  On some classes of polynomials
orthogonal on arcs of the unit circle connected with symmetric
orthogonal polynomials on an interval}, J.Approx.Theory, {\bf 94},
(1998), 73--106.

\bi{Zhe2} A.Zhedanov, {\it On the Polynomials Orthogonal on
Regular Polygons}, J.Approx.Theory, {\bf 97} (1999), 1--14.

\bi{Zh_LBP} A.Zhedanov, {\it The "classical" Laurent biorthogonal
polynomials} J. Comput. Appl. Math., {\bf 98} (1998), 121--147.

 \eb

\end{document}